\documentclass[12pt]{amsart}
\usepackage{amsmath,amscd,amssymb}

\begin{document}


\theoremstyle{plain} 
\newtheorem{thm}{Theorem}[section]
\newtheorem{cor}[thm]{Corollary}
\newtheorem{lem}[thm]{Lemma}
\newtheorem{prop}[thm]{Proposition}

\theoremstyle{definition}
\newtheorem{ex}[thm]{Example}
\newtheorem{rem}[thm]{Remark}
\newtheorem{defn}[thm]{Definition}
\newtheorem{ques}[thm]{Question}
\newtheorem{quess}[thm]{Questions}


\def \dl{^*}
\def \bu{^{\bullet}}
\def \xto{\xrightarrow}
\def \into{\hookrightarrow}
\def \onto{\twoheadrightarrow}

\def \gar{\mathbb{G}_{a(r)}}
\def \ga{\mathbb{G}_a}
\def \gao{\mathbb{G}_{a(1)}}
\def \ers{{\mathcal{E}}_{r,s}}
\def \eos{{\mathcal{E}}_{1,s}}
\def \zz{\mathbb{Z}/p}

\def \bp {\beta{\mathcal{P}}}
\def \pp {{\mathcal{P}}}

\def \F{\mathbf{F}_p}

\def \al{\alpha}
\def \la{\lambda}
\def \gl{{\mathfrak{g}}}
\def \gla{\mathfrak{g_a}}

\def\ann{\operatorname{ann}}
\def\Ext{\operatorname{Ext}}
\def\End{\operatorname{End}}
\def\Lie{\operatorname{Lie}}
\def\of{\lower3pt\hbox{${}^{\circ}$}}
\def\Hom{\operatorname{Hom}}
\def\krdim{\operatorname{Kr.dim}}
\def\spec{\operatorname{Spec}}

\title[Cohomology and projectivity]{Cohomology and projectivity of modules for
finite group schemes}
\author{Christopher P. Bendel}
\address{Department of Mathematics, Statistics, and Computer Science\\
University of Wisconsin-Stout\\
Menomonie,  WI  54751\\
U.S.A.}
\email{bendelc@uwstout.edu}
\date{January 2000}
\maketitle


\section{Introduction}

Let $G$ be a finite group scheme over a field $k$, that is, an affine group 
scheme whose coordinate ring $k[G]$ is finite dimensional.  The dual algebra 
$k[G]^* \equiv \Hom_k(k[G],k)$ is then a finite dimensional cocommutative 
Hopf algebra.  Indeed, there is an equivalence of categories between finite 
group schemes and finite dimensional cocommutative Hopf algebras
(cf.~\cite{Jan1}). Further the representation theory of $G$ is equivalent to 
that of $k[G]^*$. Many familiar objects can be considered in this context. 
For example, any finite group $G$ can be considered as a finite group scheme.
In this case, the algebra $k[G]^*$ is simply the group algebra $kG$.  Over a 
field of characteristic $p > 0$, the restricted enveloping algebra $u(\gl)$ 
of a $p$-restricted Lie algebra $\gl$ is a finite dimensional
cocommutative Hopf algebra.  Also, the mod-$p$ Steenrod algebra is graded 
cocommutative so that some finite dimensional Hopf subalgebras are such 
algebras.

Over the past thirty years, there has been extensive study of the modular 
representation theory (i.e., over a field of positive characteristic $p > 0$) 
of such algebras, particularly in regards to understanding cohomology and 
determining projectivity of modules. This paper is primarily interested in
the following two questions:

\begin{quess}\label{f1} Let $G$ be a finite group scheme $G$ over a field
$k$ of characteristic $p > 0$, and let $M$ be a rational $G$-module.
\begin{itemize}
\item[(a)] Does there exist a family of subgroup schemes of $G$ which 
detects whether $M$ is projective?
\item[(b)] Does there exist a family of subgroup schemes of $G$ which 
detects whether a cohomology class $z \in \Ext_G^n(M,M)$ (for $M$ 
finite dimensional) is nilpotent?
\end{itemize}
\end{quess}  

\noindent
It is shown here that when the connected component of $G$ is unipotent 
there is a family of subgroup schemes (with simple structure) that provides 
an affirmative answer 
to both questions. These are referred to as {\em elementary} group schemes.

\begin{defn}\label{d1} Given a field $k$ of characteristic $p > 0$ and a 
pair of non-negative integers $r$ and $s$, define the {\em elementary} 
group scheme $\ers$ to be the product $\gar\times E_s$ over $k$, where
$\gar$ denotes the $r$th Frobenius kernel of the additive group scheme $\ga$,
and $E_s$ is an elementary abelian $p$-group of rank $s$ 
(considered as a finite group scheme). The 
groups $\mathbb{G}_{a(0)}$ and $E_0$ are 
identified with the trivial group. 
\end{defn}

\noindent
The precise statements of the main results are Theorems 
\ref{restriction} and \ref{proj}.  Some 
related detection results are deduced in Sections 7 and 8 as well.

This work has its roots in work on finite groups by 
J.-P.~Serre \cite{Se}, D.~Quillen \cite{Q}, L.~Chouinard \cite{Ch}, 
J.~Alperin and L.~Evens \cite{AE}, G.~Avrunin and L.~Scott \cite{AS}, and 
J.~Carlson \cite{Ca} among others.  For a finite group, the family of
elementary abelian subgroups detects projectivity and nilpotence. Specifically,
we remind the reader of two fundamental theorems which provide a model for
other results.

\begin{thm} [Chouinard \cite{Ch}] \label{Chouinard} Let $k$ be a field of 
characteristic $p > 0$, $G$ be a finite group, and $M$ be a $kG$-module. 
Then $M$ is projective as a $kG$-module if and only if it is projective 
as a $kE$-module for every elementary abelian subgroup $E \subset G$.
\end{thm}

\begin{thm} [Quillen \cite{Q}, Carlson \cite{Ca}] \label{Quillen-Carlson} 
Let $k$ be a field of characteristic $p > 0$, $G$ be a finite group, and $M$ 
be a finite dimensional $kG$-module. Then an element $z \in \Ext_{kG}^*(M,M)$ 
is nilpotent if and only if the image of $z$ under the the natural 
restriction map $\Ext_{kG}^*(M,M) \to \Ext_{kE}^*(M,M)$ is nilpotent for every
elementary abelian subgroup $E \subset G$.
\end{thm}

More recent forerunners of this work are the analogous results which 
have been obtained for restricted Lie algebras or more generally for 
infinitesimal group schemes.  An infinitesimal group scheme $H$ is one 
for which the coordinate algebra $k[H]$ is
finite dimensional and local.  One of the most important examples of such
is the set of Frobenius kernels $\{G_{(r)}\}$ of an affine group scheme
$G$, where $G_{(r)}$ is the kernel of the $r$th iterate of the Frobenius
map on $G$.  Any restricted Lie algebra $\gl$ corresponds to a certain
group scheme $G_{(1)}$, the first Frobenius kernel of some affine group
scheme (cf.~\cite{Jan1}).  In work of E.~Friedlander, A.~Sulin, and
the author \cite{SFB2}, \cite{Bend}, it was shown that the family of
subgroup schemes $\{\gar\}$ plays the role of elementary abelian subgroups
in detecting projectivity and nilpotence of cohomology for infinitesimal
group schemes.  Specifically, Theorems 2.5 and 4.1 of \cite{SFB2} are 
analogues of Theorem \ref{Quillen-Carlson} for detecting nilpotent cohomology
classes, and Proposition 7.6 of \cite{SFB2} is an analogue of Chouinard's
theorem (Theorem \ref{Chouinard}) about detecting projectivity.  However, 
the latter result holds only for finite dimensional modules, whereas 
Chouinard's theorem holds for arbitrary modules.  In \cite{Bend}, ideas 
based on
Chouinard's original arguments were used to show that the family of
subgroup schemes
$\{\gar\}$ does detect projectivity of arbitrary modules for infinitesimal
{\em unipotent} group schemes.  It remains an open problem to show this
for arbitrary infinitesimal group schemes.

Given the results for finite groups and infinitesimal group schemes,
one is naturally led to ask the questions in \ref{f1}.  In other words, is 
there a ``nice'' general family of subgroup schemes which detects projectivity
and nilpotence?

The search for an answer begins by noticing a similarity between
the fundamentally different notions of an elmentary abelian $p$-group
$E_r$ of rank $r$ and the Frobenius kernel $\gar$.  Indeed, their 
corresponding Hopf algebras  $k[E_r]^* = kE_r$ and $k[\gar]^*$ are both 
isomorphic as algebras to a truncated polynomial algebra 
$A_r = k[x_1,x_2,\dots,x_r]/(x_1^p,x_2^p,\dots,x_r^p)$ (although 
their coalgebra structures differ).  This suggests considering a general
family of group schemes whose coordinate algebras have such structure,
and leads to the above definition (\ref{d1}) of elementary group schemes. 

It is shown here that for finite group schemes whose (infinitesimal) 
connected component at the identity, $G_0$, is {\em unipotent}
the family of elementary group schemes $\{\ers\}_{r\geq 0,s\geq 0}$ plays 
the general role of elementary abelian groups.  More precisely, the first 
main result of this paper, Theorem \ref{restriction}, is a generalization of 
Theorem \ref{Quillen-Carlson} above and Theorem 2.5 of \cite{SFB2} on detecting
nilpotent cohomology classes. And the second main result, Theorem \ref{proj}, 
is a generalization of Theorem \ref{Chouinard} above and the Theorem of 
\cite{Bend} on detecting projectivity of modules.

A group scheme $G$ is unipotent if it admits an embedding as a closed 
subgroup of $U_n$ for some $n$, where $U_n$ is the subgroup of invertible 
strictly upper triangular matrices in $GL_n$ (the group scheme of all 
invertible $n\times n$ matrices).  Unipotent groups are generalizations of 
$p$-groups in that they admit only a single simple module, the trivial module 
$k$.  (Indeed, that is sometimes taken as the definition of a unipotent group 
scheme.  And if the reader prefers a Hopf algebra perspective, then unipotent
means that we are considering Hopf algebras which admit a single simple 
module.) 

The connected component of $G$ comes into play via a decomposition of a 
finite group scheme into an infinitesimal part and a finite group part.  
More precisely, if $G$ is a finite group scheme over $k$ and all of its 
points are $k$-rational, then $G$ may be identified as a
semi-direct product $G = G_0\rtimes\pi$, where $\pi = G(k)$ is the finite 
group of $k$-points of $G$ which acts on $G_0$ via group scheme 
automorphisms (cf.~\cite{FS}, \cite{Wat}).  In the proofs of the main results,
a certain injection in cohomology (Lemma \ref{injection}) leads
to an easy reduction to the case that $\pi$ is a $p$-group.  When $\pi$ is a 
$p$-group and $G_0$ is unipotent, the group $G$ is unipotent. 

Unipotent groups have other nice properties which will also be
reviewed in Section 5.  Prior to that, some necessary cohomological 
properties of elementary group schemes are developed in Sections
2 and 3, with the key result being Proposition \ref{Serre3}.  This
is a characterization of the cohomology of elementary group schemes that
generalizes Serre's characterization of the cohomology of elementary abelian
groups (cf.~Proposition \ref{Serre1}). The proofs of the main
results, rely on homomorphisms of the form $G \onto \gao$ and $G \onto \zz$.  
Some cohomological properties related to such homomorphisms are also 
developed in Sections 4 and 5.

Theorem \ref{restriction} (whose statement and proof make up the content of 
Section 6) has several almost immediate consequences which are presented in 
Section 7.  One of these, Corollary \ref{proj-cor}, is a generalization of 
Theorem \ref{Chouinard} about detecting projectivity of {\em finite 
dimensional} modules.  Indeed, for finite dimensional modules, Chouinard's 
theorem (Theorem \ref{Chouinard}) follows via a general argument
from Theorem \ref{Quillen-Carlson}.  The key fact is that for a finite 
dimensional cocommutative Hopf algebra $A$ over $k$ and a finite dimensional 
$A$-module $M$, $\Ext_{A}^*(M,M)$ is finitely generated as a module over (the 
Noetherian ring) $H^*(A,k) = \Ext_{A}^*(k,k)$ (cf.~\cite{FS}). The bulk of
Section 8 consists of the proof of Theorem \ref{proj} which gives the
detection of projectivity for arbitrary modules.  In the final result of the 
paper (Corollary \ref{xp-proj}), Dade's lemma \cite{Dade} is used to derive 
an alternate detection result.

Having outlined the contents of the paper, we note some additional motivation
for answering the questions in \ref{f1}.  First, 
Theorem \ref{Quillen-Carlson} is a 
key component in identifying the 
{\em cohomological support varieties} of finite groups.  For any finite 
dimensional cocommutative Hopf algebra $A$ over a field $k$, the cohomology 
ring $H^*(A,k) = \Ext_A^*(k,k)$ is graded commutative so that the even 
dimensional portion $H^{2*}(A,k)$ is a commutative ring.  The cohomological
variety (or properly scheme) of $A$ is defined to be the prime ideal spectrum
$\spec H^{2*}(A,k)$.  Further, for any finite dimensional $A$-module $M$, 
consider the kernel $J_A(M)$ of the map
$$
\Ext_A^{2*}(k,k) \xto{\otimes M} \Ext_A^*(M,M).
$$
This is a homogeneous ideal which defines a homogeneous subscheme of 
$\spec H^{2*}(A,k)$ called the support variety of $M$.  These 
definitions were originally made in terms of the maximal ideal spectrum 
(over algebraically closed fields) in which case they are honest
varieties.  Within the literature, it has been noted that it is often more 
convenient to work instead with prime ideals, while maintaining the 
traditional terminology. 

This study of cohomology via varieties was extended to restricted Lie 
algebras by E.~Friedlander and B.~Parshall \cite{FP1}, \cite{FP2}, as well as 
by J.~Jantzen \cite{Jan2}, and later to arbitrary infinitesimal group 
schemes by E.~Friedlander, A.~Suslin, and the author \cite{SFB1, SFB2}.
In that work, the detection properties of the family of subgroup
schemes $\{\gar\}$ was crucial for identifying the varieties.  This suggests
that a substantial theory of cohomological varieties for arbitrary 
finite group schemes requires a family of subgroup schemes (with simple
cohomological structure) which detects nilpotence.

Further, D.~Benson, J.~Carlson, and J.~Rickard \cite{BCR1, BCR2} have 
recently developed a theory of varieties for arbitrary  modules 
(i.e., including infinite dimensional modules) which seems to be
fruitful for studying finite dimensional modules.  
A family of subgroup schemes which detects both nilpotence in cohomology and 
projectivity of arbitrary modules is important for an attempt to extend 
these ideas to more general finite group schemes or Hopf algebras, because 
some of their results depend on knowing that Theorem \ref{Chouinard} holds 
for arbitrary modules.  Indeed, the recent work of 
M.~Hovey and J.~Palmieri \cite{HP} extends the definitions of these varieties
to some more general finite dimensional graded cocommutative Hopf algebras. 

While the work here is from a group scheme 
perspective, these results could also be 
stated in terms of finite dimensional cocommutative Hopf algebras.  From the 
latter perspective, the results say that nilpotence in cohomology and 
projectivity of modules can be detected by certain {\em elementary} Hopf 
subalgebras, that is, Hopf subalgebras which have the algebra structure of a
truncated polynomial algebra $k[x_1,x_2,\dots,x_n]/(x_1^p,x_2^p,\dots,x_n^p)$.
We leave the precise formulations to the interested reader, but point out 
that others have studied these questions from a Hopf algebra perspective, in 
fact, from the more general perspective of graded Hopf algebras. 

A few examples are work of C.~Wilkerson \cite{Wil}, J.~Palmieri \cite{Pal}, 
D.~Nakano and J.~Palmieri \cite{NP}, and   M.~Hovey and J.~Palmieri \cite{HP}.
The motivation for considering graded algebras comes from topological 
interests in the Steenrod algebra which is graded.  While any Hopf algebra 
can be trivially graded (by concentrating the algebra in degree 0), a 
generic graded Hopf algebra may only be {\em graded}
cocommutative and not honestly cocommutative.  Such algebras do not 
correspond to finite group schemes and are not in our realm of consideration. 
Even if one ignores the grading, one must be careful about issues of whether 
modules and subalgebras are graded or not.  For example,  working over graded 
Hopf algebras, C.~Wilkerson \cite{Wil} showed that elementary Hopf 
subalgebras did not necessarily detect nilpotence in the
cohomology ring of the restricted enveloping algebra of a restricted Lie
algebra.  However, his example does not contradict our results because his
example is for a graded Lie algebra and he  
only considers graded subalgebras, while Theorem \ref{restriction} would 
require looking at a larger family of elementary subalgebras.

Wilkerson's work led to the more general notion of {\em quasi-elementary} 
Hopf algebras. This family of Hopf algebras was used successfully by 
J.~Palmieri \cite{Pal} to obtain an analogue of Theorem \ref{Quillen-Carlson} 
and many other results in work with D.~Nakano \cite{NP} and M.~Hovey \cite{HP}
for graded (connected) finite dimensional cocommutative Hopf
algebras.  While some of these results would apply to our situation, they 
remain somewhat unsatisfying because, while quasi-elementary Hopf algebras 
can be defined cohomologically, they have only been concretely described in 
some cases.  For a group algebra, they simply correspond to
elementary abelian groups.  For a finite dimensional Hopf subalgebra of the 
mod-2 Steenrod algebra, which is honestly cocommutative, the quasi-elementary 
subalgebras are in fact elementary. In that case, Theorem \ref{restriction}
is consistent with \cite{Pal} in claiming that
nilpotence can be detected by elementary Hopf subalgebras.

{\bf Acknowledgments.} The author gratefully acknowledges the support of the 
University of Notre Dame where this work began, thanks John Palmieri for the 
discussions about Hopf algebras and varieties which partially inspired this 
work, thanks Bill Dwyer for the suggestions on the proof of Lemma 
\ref{injection}, thanks Eric Friedlander for the comments
and suggestions on earlier versions of this work, and thanks the referee
for the helpful comments.


\section{Cohomology of Elementary Group Schemes}

For the remainder of the paper, $k$ will denote a field of characteristic 
$p > 0$ and group schemes will be defined over $k$ (unless otherwise 
specified). As group schemes of the form $\ers = \gar\times E_s$ will play a 
prominent role in the sequel, we make note of the structure of their
cohomology rings.  Of course we have 
$$H^*(\ers,k) = H^*(\gar,k)\otimes H^*(E_s,k).$$
And the cohomology of both $\gar$ and $E_s$ are well known, the first via 
\cite{CPSvdK} and the second is a classical result (see for example 
\cite{Ben} or \cite{Ev}). Indeed, their algebra
structures are identical when $r = s$, yet it will be important to 
distinguish the two pieces, and so we use different notation.  

\begin{prop}\label{ers} If $p \neq 2$, then $H^*(\gar,k) = k[x_1,x_2,\dots,
x_r]\otimes\Lambda(\lambda_1,\lambda_2,\dots,\lambda_r)$ and $H^*(E_s,k) = 
k[z_1,z_2,\dots,z_s]\otimes\Lambda(y_1,y_2,\dots,y_s)$, where the $x_i$ and 
$z_i$ have degree 2 and the $\lambda_i$ and $y_i$ have degree 1.  If $p = 2$, 
then we simply have $H^*(\gar,k) = k[\lambda_1,\lambda_2,\dots,\lambda_r]$ 
and $H^*(E_s,k) = k[y_1,y_2,\dots,y_s]$ with $\lambda_i$ and $y_i$ again 
having degree 1.  In the latter case, we set $x_i = \lambda_i^2$ and 
$z_i = y_i^2$ for each $i$.
\end{prop}

The cohomology ring of any finite-dimensional cocommutative Hopf algebra 
admits the action of Steenrod operations $\bp^i$ and $\pp^i$ (or $Sq^i$) 
(cf.~\cite{May}).  Of essential importance for our proofs will be the action 
of the Steenrod operations on $H^*(\ers,k)$.  We
summarize in the following Lemma the actions on the above generators, and 
note that the actions differ somewhat on the two pieces.  The action on the 
$E_s$ generators is well known (cf.~\cite{Ben}, \cite{May}).  For $\gar$ we 
refer the reader to \cite[Prop. 1.7]{SFB2}. 
Indeed, in the case of finite groups, the operation $\bp^i$ is simply the 
composition of $\pp^i$ with the Bockstein homomorphism $\beta$.

\begin{lem}\label{steenrod} Let $r$ and $s$ be non-negative integers. The 
action of the Steenrod operations on the generators of $H^*(\ers,k)$ as 
given in Proposition \ref{ers} are as follows:
\begin{itemize}
\item For $p > 2$,
\begin{itemize}
\item[(a)] $\pp^j(x_i) = 0 = \pp^j(z_i)$, $j > 1$; 
$\pp^j(\lambda_i) = 0 = \pp^j(y_i)$, $j \geq 1$;
\item[]$\bp^j(x_i) = 0 = \bp^j(z_i)$, $j \geq 1$;
$\bp^j(\lambda_i) = 0 = \bp^j(y_i)$, $j \geq 1$.
\item[(b)] $\pp^0(x_i) = x_{i+1}$, $\pp^0(z_i) = z_i$, 
$\pp^0(\lambda_i) = \lambda_{i+1}$, $\pp^0(y_i) = y_i$.
\item[(c)] $\pp^1(x_i) = x_i^p$, $\pp^1(z_i) = z_i^p$.
\item[(d)] $\bp^0(x_i) = 0 = \bp^0(z_i)$, $\bp^0(\lambda_i) = -x_i$, 
$\bp^0(y_i) = z_i$.
\item[(e)] $\pp^m(x_r^{\ell}) = x_r^{p\ell}$ if $m = \ell$; 
$\pp^m(x_r^{\ell}) = 0$ if $m \neq \ell$.
\end{itemize}
\item For $p = 2$,
\begin{itemize}
\item[(a)] $Sq^j(x_i) = 0 = Sq^j(z_i)$, $j > 2$; 
$Sq^j(\lambda_i) = 0 = Sq^j(y_i)$, $j > 1$.
\item[(b)] $Sq^0(x_i) = x_{i+1}$, $Sq^0(z_i) = z_i$, 
$Sq^0(\lambda_i) = \lambda_{i+1}$, $Sq^0(y_i) = y_i$.
\item[(c)] $Sq^1(x_i) = 0 = Sq^1(z_i)$, $Sq^1(\lambda_i) = x_i$, 
$Sq^1(y_i) = z_i$.
\item[(d)] $Sq^2(x_i) = x_i^2$, $Sq^2(z_i) = z_i^2$.
\item[(e)] $Sq^m(x_r^{\ell}) = x_r^{2\ell}$ if $m = 2\ell$; 
$Sq^m(x_r^\ell) = 0$ if $m \neq 2\ell$.
\end{itemize}
\end{itemize}
\end{lem}


\section{Serre's Theorem}

An essential ingredient in proving both Theorems \ref{Chouinard} and
\ref{Quillen-Carlson} is Serre's cohomological characterization of 
elementary abelian $p$-groups:

\begin{thm} [Serre \cite{Se}] \label{Serre} Let $k$ be a field of 
characteristic $p > 0$ and $G$ be a finite $p$-group.  If $G$ is not 
elementary abelian, then there exists a finite family of non-zero elements 
$v_1,v_2,\dots,v_m \in H^1(G,k)$ such that the element 
$u = \prod\beta(v_i) \in H^2(G,k)$ is zero, where $\beta$ denotes the 
Bockstein homomorphism.
\end{thm}

For our purposes it is more useful to work with 
a ``dual'' statement about the cohomology of elementary abelian groups.  
We present this
in Proposition \ref{Serre1} followed by an analogue and a generalization.
The proof of Proposition \ref{Serre1} is essentially contained in 
Serre's original proof of Theorem \ref{Serre} (which was
stated for the finite field $k = \F$).  Throughout this section we use the 
notation of 
Proposition \ref{ers}.  The following fact is a special 
case of the Corollaire to Proposition (1) of \cite{Se}.

\begin{prop}[Serre \cite{Se}]\label{fp} Let $s$ be a positive integer and
$I \subset \F[z_1,z_2,\dots,z_s] \subset H^*(E_s,\F)$ 
be a non-zero homogeneous ideal in $\F[z_1,z_2,\dots,z_s]$ which is stable 
under the Steenrod operations.  Then there exists a finite family $\{u_i\}$ 
of elements in $\F[z_1,z_2,\dots,z_s]$, each of which is a non-zero linear 
combination of the $\{z_j\}$, such that the product 
$\prod u_i \in \F[z_1,z_2,\dots,z_s]$ lies in $I$.
\end{prop}

\begin{proof} Let $A = \F[z_1,z_2,\dots,z_s]$ and consider the operation 
$\theta: A \to A$ defined
by $\theta = \sum_{i = 0}^{\infty}\pp^i$ (or the sum of the $Sq^i$ in
the case $p = 2$).  It follows from the Cartan formula,
$\pp^{\ell}(u\cdot v) = \sum_{n = 0}^{\ell}\pp^n(u)\cdot\pp^{\ell-n}(v)$,
that $\theta$ is an algebra homomorphism.  Moreover, from Lemma 
\ref{steenrod}, it is the algebra homomorphism defined by
$\theta(z_i) = z_i + z_i^p$ for each $i$.  Finally, since $I$ is stable
under the Steenrod operations, it is stable under $\theta$, and so
the claim follows from the Corollaire to Proposition (1) of \cite{Se}.
\end{proof}

The following lemma allows us to extend from the finite field $\F$ to
an arbitrary field $k$ of characteristic $p > 0$.

\begin{lem}\label{ideal} Let $k$ be a field of characteristic $p > 0$
and $s$ be a positive integer.
If $I \subset k[z_1,z_2,\dots,z_s] \subset H^*(E_s,k)$ is a non-zero 
homogeneous ideal in $k[z_1,z_2,\dots,z_s]$ which is stable under
the Steenrod operations, then 
$I' = I \cap \F[z_1,z_2,\dots,z_s]$ is also non-zero, where 
the algebra $\F[z_1,z_2,\dots,z_s]$ is naturally embedded in the
algebra $k[z_1,z_2,\dots,z_s]$.
\end{lem}

\begin{proof} Let
$$
f = \sum_{i = 1}^{n}a_iz_1^{j_{i,1}}z_2^{j_{i,2}}\cdots z_s^{j_{i,s}}
$$
be a non-zero homogeneous polynomial in $I$.  Choose such an 
$f$ that has lowest possible degree and for which the number
of summands, $n$, is least.  Having made such a choice, it suffices
to assume that $a_1 = 1$.
From Lemma \ref{steenrod} and the Cartan formula 
($\pp^0(u\cdot v) = \pp^0(u)\cdot\pp^0(u)$), we see that
$$
\pp^0(f) = \sum_{i = 1}^{n}a_i^pz_1^{j_{i,1}}z_2^{j_{i,2}}\cdots z_s^{j_{i,s}}.
$$
Since $a_1 = 1$, the difference 
$$
\pp^0(f) - f = \sum_{i = 2}^{n}(a_i^p - a_i)z_1^{j_{i,1}}z_2^{j_{i,2}}\cdots z_s^{j_{i,s}}
$$
has fewer summands than $f$.  By assumption $\pp^0(f)$ and hence 
$\pp^0(f) - f$ also lie in $I$. By the minimality of $f$, $\pp^0(f) - f$
must be zero.  That is, we must have $a_i^p = a_i$ for each $1 \leq i \leq n$.
In other words, $f$ lies in the subalgebra $\F[z_1,z_2,\dots,z_s]$ and so
the ideal $I'$ is non-zero.
\end{proof}

\begin{prop}[Serre \cite{Se}]\label{Serre1} Let $k$ be a field of 
characteristic $p > 0$, $s$ be a positive integer,  and 
$I \subset H^*(E_s,k)$ be a homogeneous ideal which is stable under the 
Steenrod operations.  If $I$ contains a non-zero element of
degree two, then there exists a finite family $\{u_i\}$ of elements in 
$H^2(E_s,k)$, each of which is a non-zero linear combination of the 
$\{z_j\}$, such that the product $\prod u_i \in H^*(E_s,k)$ lies in $I$.
\end{prop}

\begin{proof} We have 
$H^*(E_s,k) = k[z_1,z_2,\dots,z_s]\otimes\Lambda(y_1,y_2,\dots,y_s)$. 
Consider the set 
$I' = I \cap k[z_1,z_2,\dots,z_s] \subset k[z_1,z_2,\dots,z_s]$.  
This is evidently an ideal in $k[z_1,z_2,\dots,z_s]$ and moreover 
homogeneous since $I$ was.  If $I'$ is non-zero, by Lemma \ref{ideal},
the ideal $I'' = I \cap \F[z_1,z_2,\dots,z_s] \subset \F[z_1,z_2,\dots,z_s]$ 
is also non-zero as well as homogeneous and Steenrod stable.  Applying 
Proposition \ref{fp} to $I''$, the result follows.

Let $z = \sum_{1 \leq i < j \leq s}a_{i,j}y_iy_j + 
\sum_{1 \leq \ell \leq s}b_{\ell}z_{\ell}$
with $a_{i,j}, b_{\ell} \in k$ be a non-zero degree two element in $I$.  
If each $a_{i,j} = 0$, then $z$ is an element of the desired form and we're 
done.  If not, as in 
the proof of Proposition (4) of \cite{Se}, apply the operation $\bp^1\bp^0$ 
(or $Sq^1Sq^2Sq^1$ in the case $p = 2$) to $z$.
From the Cartan formula and Lemma \ref{steenrod}, we have 
$\bp^1\bp^0(z) = \sum_{i < j}a_{i,j}^{p^2}(z_i^pz_j - z_iz_j^p)$, which is 
non-zero.  Clearly this lies in $k[z_1,z_2,\dots,z_s]$ and by the
stability assumption it lies in $I$ and hence in $I'$.
\end{proof}

In the case of infinitesimal group schemes, there is a similar sort of 
result which was used in \cite{SFB2}.  For the reader's convenience, we 
restate it here.

\begin{prop} [Proposition 1.7 of \cite{SFB2}]\label{Serre2} Let $k$ be a 
field of characteristic $p > 0$, $r$ be a positive integer, and 
$I \subset H^*(\gar,k)$ be a $k$-submodule stable with respect to the 
Steenrod operations.  If $I$ contains a non-zero
element of degree two, then some power of $x_r$ lies in $I$. 
\end{prop}

We now combine these into a statement about $\ers = \gar\times E_s$.

\begin{prop} \label{Serre3} Let $k$ be a field of 
characteristic $p > 0$, $r$ and $s$ be non-negative integers, and 
$I \subset H^*(\ers,k)$ be a homogeneous ideal which is stable under the 
Steenrod operations.  If $I$ contains a non-zero element of degree
two, then there exists an integer $m$ and a finite family $\{u_i\}$ of 
elements in $H^2(E_s,k) \subset H^2(\ers,k)$, each of which is a non-zero 
linear combination of the $\{z_j\}$, such that the product 
$x_r^m\prod u_i \in H^*(\ers,k)$ lies in $I$.
\end{prop}

\begin{proof} Consider the ideals $I_1 = I \cap H^*(E_s,k)$ and 
$I_2 = I \cap H^*(\gar,k)$. These ideals satisfy the structural hypotheses 
of Propositions \ref{Serre1} and \ref{Serre2} respectively.  By assumption 
$I$ contains some non-zero element $u$ of degree 2.  Generically,
$u$ has the following form:
$$
u = \sum_{1 \leq i < j \leq r}a_{i,j}\lambda_i\lambda_j + 
\sum_{1 \leq j \leq r}b_jx_j +
\sum_{1 \leq i \leq r,1 \leq j \leq s}c_{i,j}\lambda_iy_j +
\sum_{1 \leq i < j \leq s}d_{i,j}y_iy_j + \sum_{1 \leq j \leq s}e_jz_j
$$
for some constants $a_{i,j}, b_j, c_{i,j}, d_{i,j}, e_j \in k$ which are not 
all zero.

We remind the reader of Lemma \ref{steenrod} which lists the actions of the 
Steenrod operations on the generators. We consider the case that $p > 2$ and 
simply make note of the appropriate operations in the case $p = 2$.  The 
interested reader may fill in the details of the latter case.  From Lemma 
\ref{steenrod}(b) and the Cartan formula:
$$
\pp^{\ell}(v\cdot w) = \sum_{n = 0}^{\ell} \pp^n(v)\cdot\pp^{\ell - n}(w),
$$ 
we first observe that repeated application of $\pp^0$ (resp.~$Sq^0$) to $u$
results in an element of the form
$u' = \sum_{i < j}d_{i,j}^{p^{\ell}}y_iy_j + \sum_{j}e_j^{p^{\ell}}z_j$
since the other terms are eventually killed by  $\pp^0$ (resp.~$Sq^0$). 
By the assumption that $I$ is stable under the action of the Steenrod 
operations, $u'$ lies in $I$ and hence in $I_1$. So if at least one of the 
$d_{i,j}$ or $e_j$ is non-zero, applying Proposition
\ref{Serre1} we're done.

Hence we may assume that $u$ has the form
$$
u = \sum_{i < j}a_{i,j}\lambda_i\lambda_j + \sum_{j}b_jx_j + 
\sum_{i,j}c_{i,j}\lambda_iy_j
$$
with not all coefficients being zero. As already noted, repeated 
application of $\pp^0$ (resp.~$Sq^0$) will eventually kill such a $u$.  
Stopping at the last point before we get zero,
we may further assume that $u$ has the form
$$
u = \sum_{i < r}a_{i,r}\la_i\la_r + b_rx_r + \sum_{j}c_{r,j}\lambda_ry_j
$$
with not all coefficients being zero. If every $c_{r,j} = 0$, then 
$u$ lies in $I_2$ and we are done by Proposition \ref{Serre2}.  Hence we 
may assume that some $c_{r,j}$ is non-zero.  Applying
$\bp^0$ (resp.~$Sq^1$) and using the Cartan formula:
$$
\bp^{\ell}(v\cdot w) = \sum_{n = 0}^{\ell}\bp^{n}(v)\cdot\pp^{\ell - n}(w) 
+ (-1)^{\dim v}\pp^{n}(v)\bp^{\ell - n}(w),
$$
along with Lemma \ref{steenrod}(b,d), we get
$$
\bp^0(u) = \sum_{i < r}a_{i,r}^p\lambda_{i+1}x_r - 
\sum_{j}c_{r,j}^px_ry_j \in I.
$$ 
And further applying $\bp^1$ (resp.~$Sq^3$) to this, using Lemma 
\ref{steenrod}(a,b,c,d), we get
$$
\bp^1\bp^0(u) = -  \sum_{i < r}a_{i,r}^{p^2}x_{i+1}x_r^p - 
\sum_{j}c_{r,j}^{p^2}x_r^pz_j \in I.
$$
Now we apply $\pp^p$ (resp.~$Sq^4$) to this.  With the Cartan formula and 
Lemma \ref{steenrod}(a,b,e), we get
$$
\pp^p\bp^1\bp^0(u) = -\sum_{i < r}a_{i,r}^{p^3}x_{i+2}x_r^{p^2} -
\sum_{j}c_{r,j}^{p^3}x_r^{p^2}z_j \in I.
$$
Successively applying $\pp^{p^2}$, $\pp^{p^3}$, \dots (resp.~$Sq^{2^3}$, 
$Sq^{2^4}$, \dots) we eventually conclude that $I$ contains an element of 
the form $\sum_{j}c_jx_r^{p^t}z_j = x_r^{p^t}\left(\sum_{j}c_jz_j\right)$ 
with some $c_j \neq 0$, which is of the desired form.
\end{proof}


\section{Homomorphisms}

In order to show that the collection of elementary subgroup schemes 
detect nilpotence and projectivity, we make extensive use of homomorphisms 
of the form $\phi: G \to \zz$ and $\phi: G \to \gao$.  Although the 
cohomology rings of $\zz$ and $\gao$ are identical, as in Proposition 
\ref{ers}, we will continue to use different notation to distinguish the two: 
$H^*(\zz,k) = k[z_1]\otimes \Lambda(y_1)$ (or $k[y_1]$ if $p = 2$,
with $z_1 = y_1^2$) and similarly 
$H^*(\gao,k) = k[x_1]\otimes\Lambda(\lambda_1)$ (or $k[\lambda_1]$ if 
$p = 2$, with $x_1 = \lambda_1^2$). Of particular importance will be the
image $\phi^*(z_1)$ (resp.~$\phi^*(x_1)$) of $z_1$ (resp.~$x_1$) under
the induced map in cohomology $\phi^*: H^*(\zz,k) \to H^*(G,k)$ 
(resp.~$\phi^*: H^*(\gao,k) \to H^*(G,k)$).  The significance of such 
homomorphisms is partially seen in the following key properties.

An essential ingredient in the proof of Theorem \ref{restriction} is the 
following generalization of the Quillen-Venkov Lemma \cite{QV}, whose proof 
is the same as that of \cite[Proposition 2.3]{SFB2} (see also \cite{Ben}).  
A unital rational $G$-algebra $\Lambda$ is a unital $k$-algebra which admits 
a compatible structure of a rational $G$-module.  That is, $1_{\Lambda}$ 
lies in  $\Lambda^G$ and the multiplication map 
$\Lambda\otimes_k\Lambda \to \Lambda$ is a map of rational $G$-modules.  
For example, take $\Lambda = \Hom_k(M,M)$ for a
rational $G$-module $M$.  For such an algebra, we denote by 
$\rho_{\Lambda}: k \to \Lambda$
the canonical homomorphism which defines $1_{\Lambda}$.  Further, let 
$\rho_{\Lambda}^*: H^*(G,k) \to H^*(G,\Lambda)$ denote the induced map 
in cohomology.

\begin{prop}\label{QV} Let $k$ be a field of characteristic $p > 0$, $G$ be 
a finite group scheme over $k$, $\Lambda$ be an associative, unital 
rational $G$-algebra, and $\phi: G \to \zz$ 
{\rm{(}}resp.~$\phi: G \to \gao${\rm{)}} be a non-trivial homomorphism of
group schemes over $k$.  Let $z \in H^n(G,\Lambda)$ satisfy 
$z|_{\ker \phi} = 0$.  Then $z^2$ is divisible by 
$\rho_{\Lambda}^*(\phi^*(z_1))$
{\rm{(}}resp.~$\rho_{\Lambda}^*(\phi^*(x_1))${\rm{)}} $\in H^2(G,\Lambda)$.
\end{prop}

One of the ingredients in proving Proposition \ref{QV} is the following well 
known property of the action of $H^*(H,k)$ on $H^*(H,Q)$, where $H$ denotes
either $\zz$ or $\gao$, and $Q$ is any $H$-module 
(cf.~\cite{Ch}, \cite[2.3]{SFB2}): 
the action of $z_1$ or $x_1$ (as appropriate) induces a periodicity 
isomorphism $H^j(H,Q) \xto{\sim} H^{j+2}(H,Q)$ for all $j > 0$.  
This fact is also the key ingredient in proving the following lemma which 
gives a condition for a periodicity isomorphism for $H^*(G,M)$ (for a 
rational $G$-module $M$) with respect to the action of $\phi^*(z_1)$ or
$\phi^*(x_1)$.  This is a simple generalization of Lemma 1.2 of \cite{Ch} 
(see also \cite{Bend}) and will be used several times in the proof of 
Theorem \ref{proj}.

\begin{lem}[Chouinard \cite{Ch}]\label{action} Let $k$ be a field of 
characteristic $p > 0$, $G$ be a finite group scheme over $k$, and 
$M$ be a rational $G$-module. Let $\phi: G \to \zz$
{\rm{(}}resp.~$\phi: G \to \gao${\rm{)}} be a non-trivial homomorphism of 
group schemes over $k$ and let $N$ denote the kernel of $\phi$. Suppose 
that $H^j(N,M) = 0$ for all $j > 0$, then
under the action of $H^*(G,k)$ on $H^*(G,M)$, the action of $\phi^*(z_1)$
{\rm{(}}resp.~$\phi^*(x_1)${\rm{)}} induces a periodicity isomorphism 
$H^j(G,M) \xto{\sim} H^{j + 2}(G,M)$ for all $j > 0$.
\end{lem}

\begin{proof} Let $H$ denote either $\zz$ or $\gao$ and $\phi: G \onto H$ be
the given homorphism.  Since $\phi$ is non-trivial we have a short exact 
sequence of group schemes over $k$:
$$
1 \to N \to G \to H \to 1
$$
and a corresponding Lyndon-Hochschild-Serre spectral sequence
$$
E_2^{p,q}(M) = H^p(H,H^q(N,M)) \Rightarrow H^{p+q}(G,M).
$$
By assumption, the spectral sequence collapses to
$$
E_2^{p,0}(M) = H^p(H,H^0(N,M)) \Rightarrow H^{p+q}(G,M)
$$
and hence we have an isomorphism $H^*(G,M) \simeq H^*(H,Q)$ where 
$Q = H^0(N,M)$.  The cohomology ring $H^*(H,k)$ acts on the spectral 
sequence with the action on the abutment via $\phi^*$.  Hence, the 
periodicity isomorphism follows from that of $z_1$  or $x_1$ on $H^*(H,Q)$.
\end{proof}


\section{Cohomology Facts}

The homomorphisms discussed in the previous section are of interest because 
they may be considered as elements of $H^1(G,k)$.  More precisely, for any 
affine group scheme $G/k$, we have $H^1(G,k) = \Hom_{Gr/k}(G,\ga)$ as abelian 
groups (cf.~\cite[Lemma 1.1 (a)]{SFB2}).  This can be refined in the case 
that $G$ may be identified as a semi-direct product $G = G_0\rtimes\pi$ of 
an infinitesimal group scheme $G_0$ with a finite group $\pi$
which acts on $G_0$ via group scheme automorphisms.  In such a situation,
let $\Hom_{Gr/k}(G_0,\ga)^{\pi}$ denote those homomorphisms of group schemes 
over $k$ which are preserved by the action of $\pi$.  That is to say, those 
$\phi: G_0 \to \ga$ such that $\phi(x(g)) = \phi(g)$ for all $x \in \pi$ and 
$g \in G_0$, where $x(g)$ denotes the image of $g$ under the action of $x$ 
on $G_0$. 

\begin{lem}\label{h1} Let $G$ be an affine group scheme over $k$ which may be 
identified as a semi-direct product $G = G_0\rtimes\pi$ of an infinitesimal 
group scheme $G_0$ with a finite group $\pi$.  Then 
$H^1(G,k) = \Hom_{Gr/k}(G_0,\ga)^{\pi}\times\Hom_{Gr/k}(\pi,\ga)$
{\rm{(}}as abelian groups{\rm{)}}.  
\end{lem}

\begin{proof} Any homomorphism $\phi: G \to \ga$ determines two 
homomorphisms $\phi_1: G_0 \to \ga$ by $\phi_1(g) = \phi((g,1_{\pi}))$ and 
$\phi_2: \pi \to \ga$ by $\phi_2(x) = \phi((1_{G_0},x))$.  Further, 
$\phi_1$ necessarily preserves the action of $\pi$.  Given any $g
\in G_0$ and $x \in \pi$, since $\phi$ is a homomorphism, we have 
$$
\phi((x(g),x)) = \phi((x(g),1_{\pi})\cdot(1_{G_0},x)) = \phi((x(g),1_{\pi})) 
+ \phi((1_{G_0},x)) = \phi_1(x(g)) + \phi_2(x).
$$
On the other hand, we also have
$$
\phi((x(g),x)) = \phi((1_{G_0},x)\cdot(g,1_{\pi})) = \phi((1_{G_0},x)) +
\phi((g,1_{\pi})) = \phi_2(x) + \phi_1(g).
$$
Hence, we must have $\phi_1(x(g)) = \phi_1(g)$ as claimed.

Conversely, given homomorphisms $\phi_1 \in \Hom_{Gr/k}(G_0,\ga)^{\pi}$ and 
$\phi_2 \in \Hom_{Gr/k}(\pi,\ga)$, define a map $\phi: G \to \ga$ by 
$\phi((g,x)) = \phi_1(g) + \phi_2(x)$.  For any pair 
$(g_1,x_1), (g_2,x_2) \in G$, we check that $\phi$
is in fact a homomorphism.  On the one hand, we have
\begin{align*}
\phi((g_1,x_1)\cdot(g_2,x_2)) &= \phi((g_1\cdot x_1(g_2),x_1\cdot x_2) = 
                   \phi_1(g_1\cdot x_1(g_2)) + \phi_2(x_1\cdot x_2)\\
     &= \phi_1(g_1) + \phi_1(x_1(g_2)) + \phi_2(x_1) + \phi_2(x_2).
\end{align*}
On the other hand,
$$
\phi((g_1,x_1)) + \phi((g_2,x_2)) = 
\phi_1(g_1) + \phi_2(x_1) + \phi_1(g_2) + \phi_2(x_2).
$$
Since $\phi_1(x_1(g_2)) = \phi_1(g_2)$, these agree and $\phi$ is a 
homomorphism.
\end{proof}

The first step in the proofs of the main results is to reduce to the case 
that $G$ is a unipotent group scheme.  This allows us to make use of the 
following result of \cite{SFB2} which is a generalization of a classical 
result for finite $p$-groups (cf.~\cite[Theorem 7.2.4]{Ev}).

\begin{lem}[{\cite[Lemma 1.2]{SFB2}}]\label{h1h2} Let $\psi: G \to H$ be a 
surjective homomorphism of unipotent group schemes over $k$.  If 
$\psi^*: H^1(H,k) \to H^1(G,k)$ is an isomorphism and if 
$\psi^*: H^2(H,k) \to H^2(G,k)$ is injective, then $\psi$ is an isomorphism.
\end{lem}

For unipotent group schemes there is a nice cohomological criterion 
(Proposition \ref{uniph1}) for detecting projectivity which is well known 
for $p$-groups and holds more generally. We remind the 
reader that for a finite group $G$ a $kG$-module is projective if and only 
if it is injective.  More generally, any finite dimensional cocommutative 
Hopf algebra is a Frobenius algebra (cf.~\cite{Jan1}, \cite{LS}) and
hence a module over such an algebra is in fact projective if and only if 
it is injective (cf.~\cite{FW}).  So this equivalence holds for arbitrary 
finite group schemes.  Further, the theory of finite dimensional algebras 
shows that if $k$ is the only simple module for such an
algebra (e.g., for a unipotent group scheme), then the algebra is 
indecomposable as a module over itself, and hence a module is in fact 
projective if and only if it is free.  Thus, over a finite unipotent group 
scheme, the notions of projective, injective, and free are all equivalent.

\begin{prop}[cf. \cite{Ben} or \cite{Bend}]\label{uniph1} Let $k$ be a 
field of characteristic $p > 0$ and $G$ be a finite unipotent group scheme 
over $k$.  For any rational $G$-module $M$, $M$ is projective 
{\rm{(}} = injective = free {\rm{)}} if and only if $H^1(G,M) = 0$.
\end{prop}

In order to reduce to the case of a unipotent group scheme, we use the 
injectivity of certain restriction maps in cohomology.  If $\pi$ is a finite 
group and $\pi_S \subset \pi$ is a $p$-Sylow subgroup, then for any 
$k\pi$-module $M$ the restriction map in cohomology 
$H^*(\pi,M) \to H^*(\pi_S,M)$ is an injection (cf.~\cite{Ben} or \cite{Ev}).  
This can be partially extended to the more general setting of finite group 
schemes.

\begin{lem}\label{injection} Let $k$ be a field of characteristic $p > 0$ 
and $G$ be a finite group scheme over $k$ which can be identified as a 
semi-direct product $G = G_0\rtimes\pi$ of an infinitesimal group scheme 
$G_0$ {\rm{(}}which is a closed subgroup scheme in $G${\rm{)}}
and a finite group $\pi$.  Further, let $\pi_S \subset \pi$ be a $p$-Sylow 
subgroup of $\pi$ and $M$ be a rational $G$-module.  Then the restriction 
map $H^*(G,M) \to H^*(G_S,M)$ for the embedding 
$G_S = G_0\rtimes\pi_S \into G_0\rtimes\pi = G$ is an injection.
\end{lem}

\begin{proof}  The idea is to identify $G$-cohomology groups with certain 
$\pi$-cohomology groups.  We first note that there exists a projective 
resolution $X_{\bu} \onto k$ of rational $G_0$-modules which admits a 
compatible action of the finite group $\pi$.  That is, this is also a 
complex of $\pi$-modules and for any $x \in \pi$, $g \in G_0$, and 
$m \in X_n$ (any $n$) we have $x\cdot(g\cdot m) = x(g)\cdot(x\cdot m)$.  
Consider the cobar resolution (cf.~\cite{Jan1})
$$
k \to k[G] \to k[G]^{\otimes2} \to \dots
$$
which is an injective resolution of $k$ over $G$.  Further, this is a 
sequence of rational $G$-modules and maps if $G$ is defined to act by the
left regular representation on the first $k[G]$ factor of each term 
$k[G]^{\otimes n}$ (cf.~\cite{Jan1}).   Hence $\pi \subset G$ acts
on the sequence and does so compatibly with respect to its action on $G_0$.  
Let $X_{\bu}$ be the dual complex (i.e., take the $k$-linear dual of each
module and reverse the arrows). Then $X_{\bu}$ is necessarily a  projective 
resolution of $k$ over $G$, which can be considered as a
{\em resolution} over $G_0$ on which $\pi$ acts compatibly.  Finally, since 
$G_0$ is a closed subgroup scheme of $G$, any projective $G$-module is also a 
projective $G_0$-module (cf.~\cite{Jan1}).  So $X_{\bu}$ is in fact a 
{\em projective} resolution of $k$ over $G_0$ with $\pi$ acting compatibly.

Given that such a resolution exists, the same argument as on p.~19 of 
\cite{Ev} (for a semi-direct product of finite groups) shows that there is 
an isomorphism
$$
H^*(G,M) = H^*(G_0\rtimes\pi,M) \simeq H^*(\pi,\Hom_{G_0}(X_{\bu},M))
$$
where the latter cohomology group is a {\em hypercohomology} group.  That is, 
the coefficients consist of a cochain complex of modules.  Further, this 
isomorphism is preserved under the natural embedding $\pi_S \into \pi$ so 
that the following diagram commutes:
$$
\CD
H^*(G_0\rtimes\pi,M) @>{\sim}>> H^*(\pi,\Hom_{G_0}(X_{\bu},M))\\
@V{res}VV                          @VV{res}V\\
H^*(G_0\rtimes\pi_S,M) @>{\sim}>> H^*(\pi_S,\Hom_{G_0}(X_{\bu},M)).
\endCD
$$
Finally, the left hand map will be injective if the right hand map is.
However, the injectivity of the right hand restriction map  in ordinary 
cohomology can be extended to hypercohomology. For example, one standard 
proof in ordinary cohomology (cf.~\cite{Ben} or \cite{Ev}) is
based on the fact that for a $k\pi$-module $N$, the composite
$$
H^*(\pi,N) \xto{res} H^*(\pi_S,N) \xto{Tr} H^*(\pi,N)
$$
is multiplication by the index $[\pi:\pi_S]$ which is invertible in $k$, 
where $Tr$ denotes the transfer map.  Hence the composite is an injection 
and so the restriction map is also. This same argument works with $N$ 
replaced by a cochain complex.

Alternately, by the Eckmann-Shapiro lemma, we may identify
$$
H^*(\pi_S,N) = \Ext_{\pi_S}^*(k,N) 
\simeq \Ext_{\pi}^*(k\pi\otimes_{k\pi_S}k,N),
$$ 
and the restriction map
$$
H^*(\pi,N) = \Ext_{\pi}^*(k,N) \xto{res} \Ext_{\pi_S}^*(k,N) \simeq
\Ext_{\pi}^*(k\pi\otimes_{k\pi_S}k,N)
$$
is simply the map induced from the canonical map 
$k\pi\otimes_{k\pi_S}k \to k$ which sends $x\otimes c \mapsto c$.  
This module map is split by the map $k \to k\pi\otimes_{k\pi_S}k$ which 
sends $1 \mapsto \frac{1}{|T|}\sum_{t \in T}t\otimes 1$, 
where $T \subset \pi$ is a set of left coset representatives of 
$\pi_S$ in $\pi$.  This splitting of the module map induces a splitting 
of the restriction map in cohomology when $N$ is either a module
or a complex.
\end{proof}


\section{The Restriction Theorem in Cohomology}

In this section we present a generalization of Theorem \ref{Quillen-Carlson} 
on detecting nilpotent cohomology classes, which is also an extension of 
Theorem 2.5 of \cite{SFB2} which applies to infinitesimal unipotent group 
schemes.

\begin{thm} \label{restriction} Let $k$ be a field of characteristic $p > 0$, 
$G$ be a finite group scheme over $k$, and $\Lambda$ be an associative, 
unital rational $G$-algebra {\rm{(}}as defined in Section 4{\rm{)}}. Suppose 
further that the {\rm{(}}infinitesimal{\rm{)}} connected component of the 
identity, $G_0$, of $G$ is unipotent.  If $z \in H^n(G,\Lambda)$ satisfies 
the property that for any field extension $K/k$ and any group scheme 
embedding $\nu: \ers\otimes_kK \into G\otimes_kK$ over $K$ the cohomology 
class $\nu^*(z) \in H^n(\ers\otimes_kK,\Lambda\otimes_kK)$ is nilpotent, 
then $z$ is itself nilpotent.
\end{thm}

\begin{proof} The strategy of the proof is to extend that of Theorem 2.5 
in \cite{SFB2}.  For any field extension $K/k$, since
$H^*(G\otimes_kK,\Lambda\otimes_kK) = H^*(G,\Lambda)\otimes_kK$ 
(cf.~\cite{Jan1}), if $z$ is nilpotent after base change, then it is 
necessarily nilpotent.  Hence it suffices to assume that $k$ is 
algebraically closed.  In this case, the points of $G$ are certainly
$k$-rational and so, as noted in the Introduction, we can identify $G$ as 
the semi-direct product $G_0\rtimes\pi$ where $\pi = G(k)$ is the finite 
group of $k$-points of $G$.  By assumption, $G_0$ is unipotent, but $\pi$ 
may be an arbitrary finite group.  However, by Lemma \ref{injection}, we 
may assume that $\pi$ is in fact a $p$-group, for, if $z$ is nilpotent
after restriction to the subgroup $G_0\rtimes\pi_S$ (for a $p$-Sylow 
subgroup $\pi_S$ of $\pi$), since the restriction map in cohomology is 
injective, $z$ must necessarily be nilpotent. Hence, the group $G$ may be 
assumed to be unipotent.

We now proceed by induction on $\dim_kk[G]$ and are trivially done if 
$\dim_kk[G] = 1$ or $G = \ers$, and so assume that the theorem holds for 
all groups $H$ over any field $K/k$ with $\dim_KK[H] < \dim_kk[G]$. As in 
Section 4, let $\rho_{\Lambda}^*: H^*(G,k) \to H^*(G,\Lambda)$ denote the 
map induced by the module map $\rho_{\Lambda}: k \to \Lambda$ which
defines $1_{\Lambda}$.

The next step is to reduce to the case that $\pi$ is elementary abelian.  
If not, by Serre's theorem (Theorem \ref{Serre}), there exists a finite 
product $u = \prod u_i = \prod\beta(v_i) \in H^*(\pi,k)$ which is zero in 
$H^*(\pi,k)$ for non-zero $v_i \in H^1(\pi,k)$.  Each $v_i$
can be considered as a non-zero map $\pi \to \zz$, which can be extended to 
$\phi_i: G \onto \pi \onto \zz$.  Let $N_i$ denote the kernel of $\phi_i$.  
Since $\phi_i$ is non-trivial, $\dim_kk[N_i] < \dim_kk[G]$ and so by 
induction, the restriction of $z$ to each $N_i$ is nilpotent.  Hence, by 
Proposition \ref{QV}, $z^2$ is divisible by 
$\rho_{\Lambda}^*(\phi_i^*(z_1))$ for each $i$.  Consider the canonical
projection $\psi: G \to \pi$ and the induced map
$\psi^*: H^*(\pi,k) \to H^*(G,k)$.  Since $\psi^*(u_i) = \phi_i^*(z_1)$ for
each $i$, some power of $z$ is divisible by  $\rho_{\Lambda}^*(\psi^*(u))$ 
and hence is zero.

Hence, we may assume that $\pi = E_s$ for some $s$ (possibly zero). In the 
trivial case $s = 0$, we have $G = G_0$ and the claim is precisely 
Theorem 2.5 of \cite{SFB2}.  In any case, the succeeding argument holds in 
general (with Case I being impossible in that situation).  The
remainder of the argument consists of three cases depending on
$\dim_k\Hom_{Gr/k}(G_0,\ga)^{\pi}$ and uses arguments similar to those in 
the above reduction to the case that $\pi$ is elementary abelian.

\smallskip
\noindent
{\bf CASE I:} $\dim_k\Hom_{Gr/k}(G_0,\ga)^{\pi} = 0$

Consider the natural projection $\psi: G = G_0\rtimes \pi \to \pi$ and the 
induced map $\psi^*: H^*(\pi,k) \to H^*(G,k)$.  By the assumption and 
Lemma \ref{h1}, we have an isomorphism $\psi^*: H^1(\pi,k) \to H^1(G,k)$.  
Consider the map in degree 2, $\psi^*: H^2(\pi,k) \to H^2(G,k)$.  By 
Lemma \ref{h1h2}, either $G \simeq \pi \simeq E_s$ and we are
done or $I \equiv \ker(\psi^*)$ contains some non-zero element of degree 2.  
Now, the map $\psi^*$ preserves the action of the Steenrod algebra 
(cf.~\cite{Dold}, \cite{Pr}, and also \cite{SFB2}) and hence $I$ is a 
homogeneous ideal stable under the action of the Steenrod operations.  By
Proposition \ref{Serre1}, there exists a non-zero product 
$u = \prod u_i = \prod\beta(v_i) \in H^*(\pi,k)$ which lies in $I$.  
In other words, the image of $u$ under $\psi^*$ is zero in
$H^*(G,k)$. On the other hand, arguing as above, by Proposition \ref{QV}, 
we conclude that some power of $z$ is divisible by 
$\rho_{\Lambda}^*(\psi^*(u))$ and hence is zero.

\smallskip
\noindent
{\bf CASE II:} $\dim_k\Hom_{Gr/k}(G_0,\ga)^{\pi} = 1$

Let $\theta: G_0 \to \ga$ be a representative map preserved by $\pi = E_s$.  
Since $G_0$ is infinitesimal, the image of $\theta$ lands in some $\gar$. 
If $r >1$, the composition
$$
G_0 \overset{\theta}{\to} \ga \overset{F}{\to} \ga,
$$ 
where $F$ is the Frobenius morphism, would give another linearly independent
homomorphism preserved by $\pi$.  Hence we must have $\theta: G_0 \to \gao$.  
Consider the homomorphism 
$\eta: G = G_0\rtimes E_s \xto{\theta\times Id} \gao\times E_s = \eos$ 
and the induced morphism $\eta^*: H^*(\eos,k) \to H^*(G,k)$.  Specifically, 
consider the map on $H^1$ (cf.~Lemma \ref{h1}):
\begin{align*}
H^1(\eos,k) = &\Hom_{Gr/k}(\gao,\ga)\times\Hom_{Gr/k}(E_s,\ga) \\
    &\to \Hom_{Gr/k}(G_0,\ga)^{E_s}\times\Hom_{Gr/k}(E_s,\ga) = H^1(G,k).
\end{align*}
Again, this is evidently an isomorphism and as before, either 
$G \simeq \eos$ and we're done or there exists a non-zero element of 
degree 2 in $I = \ker(\eta^*)$.  In the latter case, by
Proposition \ref{Serre3}, there exists a product 
$u = x_1^n\prod u_i \in H^*(\eos,k)$ which lies in $I$.  

Each element $u_i$ may again be identified with $\beta(v_i)$ for some 
non-zero $v_i \in H^1(E_s,k)$.  Just as above, for each $i$, we consider 
the homomorphism $\phi_i: G \onto E_s \onto \zz$ corresponding to $v_i$ and 
conclude by Proposition \ref{QV} that $z^2$ is divisible by 
$\rho_{\Lambda}^*(\eta^*(u_i))$.

On the other hand, consider the composite
$$
\phi: G = G_0 \rtimes E_s \xto{\eta} \gao\times E_s \onto \gao,
$$
where the last map is the canonical projection. By Proposition \ref{QV}, 
$z^2$ is also divisible by $\rho_{\Lambda}^*(\eta^*(x_1))$. Hence 
some power of $z$ is divisible by $\rho_{\Lambda}(\eta^*(u))$ and
thus zero since $\eta^*(u)$ is zero in $H^*(G,k)$.

\smallskip
\noindent
{\bf CASE III:} $\dim_k\Hom_{Gr/k}(G_0,\ga)^{\pi} > 1$

This case must be further divided into two, based on 
$\dim_k\Hom_{Gr/k}(G_0,\gao)^{\pi}$. As noted above, since $G_0$ is 
infinitesimal, any $\phi: G_0 \to \ga$ has image in some $\gar$.
So there is certainly one non-trivial $\pi$-preserved map $G_0 \to \gao$, 
but there need not be two linearly independent such maps.

\smallskip
{\bf CASE III(a):} $\dim_k\Hom_{Gr/k}(G_0,\gao)^{\pi} = 1$

The proof of this case is similar to the proof of Theorem 1.6 in 
\cite{SFB2}. Let $\theta: G_0 \to \gao$ be a non-trivial representative 
map preserved by $\pi$.  Since $\theta$ is preserved by the action of $\pi$, 
it can be extended to a non-trivial map $\phi: G \to \gao$ by 
$\phi((g,x)) = \theta(g)$.  The Frobenius map $F: \ga \to \ga$ induces via 
composition a map (of the same name) 
$F: \Hom_{Gr/k}(G_0,\ga) \to \Hom_{Gr/k}(G_0,\ga)$.  Further, since
any $\pi$-preserved map remains so after composition with $F$, this 
restricts to a map 
$F: \Hom_{Gr/k}(G_0,\ga)^{\pi} \to \Hom_{Gr/k}(G_0,\ga)^{\pi}$, 
and we may identify the kernel of this (restricted) map with 
$\Hom_{Gr/k}(G_0,\gao)^{\pi}$ (cf.~\cite[Lemma 1.1]{SFB2}).  By
assumption, the kernel must then be one-dimensional, and so there exists a 
non-negative integer $r$ and homomorphism 
$\zeta \in \Hom_{Gr/k}(G_0,\gar)^{\pi}$ with $F^{r-1}(\zeta) = \theta$ and 
such that the set
$\{\zeta,F(\zeta),\dots,F^{r-1}(\zeta)\}$ is a basis for the subspace
$\Hom_{Gr/k}(G_0,\ga)^{\pi} \subset H^1(G,k)$.   

Consider the homomorphism 
$\eta: G = G_0\rtimes E_s \xto{\zeta\times Id} \gar\times E_s$ and
the induced map on cohomology $\eta^*: H^*(\ers,k) \to H^*(G,k)$. 
On $H^1$ the map
\begin{align*}
H^1(\ers,k) = &\Hom_{Gr/k}(\gar,\ga)\times\Hom_{Gr/k}(E_s,\ga) \\
    &\to \Hom_{Gr/k}(G_0,\ga)^{E_s}\times\Hom_{Gr/k}(E_s,\ga) = H^1(G,k)
\end{align*}
is the identity on the right factor and on the left factor maps the basis
$\la_1, \la_2, \dots, \la_r$ to the basis 
$\zeta, F(\zeta), \dots, F^{r-1}(\zeta) = \theta$. 
Hence, this is an isomorphism and either $G \simeq \ers$ and we're done or 
there is a homomorphism 
$\eta: G  \to \ers$ for which (by Proposition \ref{Serre3}) the kernel 
of $\eta^*: H^*(\ers,k) \to H^*(G,k)$ contains an element of the form 
$u = x_r^n\prod u_i$ with $\eta^*(x_r) = \theta^*(x_1) = \phi^*(x_1)$ 
and for each $i$, $\eta^*(u_i) = \phi_i^*(z_1)$ for some non-trivial 
$\phi_i: G \to \zz$.  By Proposition \ref{QV} and the usual argument, 
some power of $z$ is divisible by $\rho_{\Lambda}^*(\eta^*(u))$ and hence 
is zero.

\smallskip
{\bf CASE III(b):} $\dim_k\Hom_{Gr/k}(G_0,\gao)^{\pi} > 1$

Let $\theta_1: G_0 \to \gao$ and $\theta_2: G_0 \to \gao$ be two linearly 
independent (and non-trivial) homomorphisms which are preserved by the 
action of $\pi$.  Since they are preserved by $\pi$, these maps can be 
extended to maps $\phi_1: G \to \gao$ by $\phi_1((g,x)) = \theta_1(g)$ and 
$\phi_2: G \to \gao$ by $\phi_2((g,x)) = \theta_2(g)$.  Clearly these remain
non-trivial and linearly independent maps $G \to \gao$.  Now, the identical 
argument as in the case $\dim_k\Hom_{Gr/k}(G,\gao) > 1$ in the proof of 
\cite[Theorem 2.5]{SFB2} may be applied with $\phi_1$ and $\phi_2$ to imply 
that $z$ is indeed nilpotent.  
\end{proof}

\begin{ques} Can this result be extended to any finite group scheme $G$?
\end{ques}

\begin{rem}  Given a homomorphism $\ers \to G$ of group schemes over $k$, 
the image of $\ers$ is necessarily a closed subgroup and moreover isomorphic 
to an elementary group scheme ${\mathcal{E}}_{r',s'}$ for some 
$r' \leq r, s' \leq s$.  Hence, Theorem \ref{restriction} (as well as 
succeeding results) could be stated either in terms of homomorphisms of 
the form $\ers \to G$ or in terms of closed subgroup schemes of the 
form $\ers$.
\end{rem}


\section{Consequences}

In this section, we note some immediate consequences of Theorem 
\ref{restriction}. First, we note a slightly weaker version of the theorem, 
which is stated in terms of Ext-groups like Theorem \ref{Quillen-Carlson}.  
Given a rational $G$-module $M$, the algebra $\Lambda \equiv \Hom_k(M,M)$ 
is an associative, unital rational $G$-algebra, and hence the theorem applies 
to $\Lambda$.  Further, there is a natural isomorphism 
$H^*(G,\Lambda) \simeq \Ext_G^*(M,M)$.  In the case that $M$ is finite 
dimensional, for any field extension $K/k$, we have 
$\Ext_{G\otimes_kK}^*(M\otimes_kK,M\otimes_kK) 
  \simeq \Ext_G^*(M,M)\otimes_kK$. 
Hence, the following is an immediate consequence of Theorem \ref{restriction}.

\begin{cor}\label{ext-version} Let $k$ be a field of characteristic $p > 0$, 
$G$ be a finite group scheme over $k$, and $M$ be a finite dimensional 
rational $G$-module. Suppose further that the {\rm{(}}infinitesimal{\rm{)}} 
connected component of the identity, $G_0$, of $G$ is unipotent.  
If $z \in \Ext_G^n(M,M)$ satisfies the property that for any field extension
$K/k$ and any group scheme embedding $\nu: \ers\otimes_kK \into G\otimes_kK$ 
over $K$ the cohomology class 
$\nu^*(z) \in \Ext^n_{\ers\otimes_kK}(M\otimes_kK,M\otimes_kK)$ is nilpotent,
then $z$ is itself nilpotent.
\end{cor}

A standard argument (cf.~\cite{Ben}, \cite{SFB2}) which is sketched below 
shows that a generalization of Chouinard's theorem (Theorem \ref{Chouinard}) 
follows from Corollary \ref{ext-version}.  

\begin{cor}\label{proj-cor} Let $k$ be a field of characteristic $p > 0$, 
$G$ be a finite group scheme over $k$, and $M$ be a finite dimensional 
rational $G$-module. Suppose further that the {\rm{(}}infinitesimal{\rm{)}} 
connected component of the identity, $G_0$, of $G$ is unipotent. Then $M$ is 
projective as a $G$-module if and only if for every field extension $K/k$ and 
closed subgroup scheme $H \subset G\otimes_kK$ with $H \simeq\ers\otimes_kK$ 
the restriction of $M\otimes_kK$ to $H$ is projective.
\end{cor}

\begin{proof} If $M$ is projective over $G$, then it remains so upon 
restriction to any closed subgroup scheme (cf.~\cite{Jan1}).  Conversely, 
given an $H \subset G\otimes_kK$ with $H \simeq \ers\otimes_kK$, if 
$M\otimes_kK$ is projective upon restriction, then
$\Ext_H^i(M\otimes_kK,M\otimes_kK) = 0$ for all $i > 0$.  Hence, by 
Corollary \ref{ext-version}, every element $z \in \Ext_G^n(M,M)$ for 
$n > 0$ is nilpotent.  Since, by \cite{FS}, $\Ext_G^*(M,M)$ is finitely 
generated over the Noetherian ring $H^{2*}(G,k) = \Ext^{2*}_G(k,k)$, 
we must have $\Ext^i_G(M,M) = 0$ for all $i > N$ for some $N$.  As the
notions of projective and injective are equivalent (see the discussion 
preceding Proposition \ref{uniph1}), it follows that $M$ must in fact be 
projective.  
\end{proof}

Let $\gl$ be a finite dimensional restricted Lie algebra over $k$ with
restricted enveloping algebra $u(\gl)$.  As mentioned in the Introduction, 
this corresponds to a certain (height 1) infinitesimal group scheme.  
Theorem \ref{restriction} applied to the corresponding group scheme says 
that nilpotence is detected by the collection of subalgebras 
$u(\langle x\rangle) \subset u(\gl\otimes_kK)$ for each
$p$-nilpotent element $x \in \gl\otimes_kK$.  Here $u(\langle x \rangle)$
denotes the subalgebra generated by $x$ and is simply isomorphic to 
$K[x]/(x^p)$ as an algebra.  Using this special case, we can obtain a 
restriction theorem in terms of subalgebras of the form $K[x]/(x^p)$.

\begin{cor}\label{xp-nilp}  Let $k$ be a field of characteristic $p > 0$, 
$G$ be a finite group scheme over $k$, and $M$ be a finite dimensional 
rational $G$-module. Suppose further that the {\rm{(}}infinitesimal{\rm{)}} 
connected component of the identity, $G_0$, of $G$ is unipotent.  If 
$z \in \Ext_G^n(M,M)$ satisfies the property that for any field extension
$K/k$ and any embedding of algebras 
$\nu: B \into K[G\otimes_kK]^*$ with $B \simeq K[x]/(x^p)$ 
the cohomology class 
$\nu^*(z) \in \Ext^n_{B}(M\otimes_kK,M\otimes_kK)$ is nilpotent,
then $z$ is itself nilpotent.
\end{cor}

\begin{proof}  For any closed subgroup scheme of the form 
$\ers\otimes_kK \subset G\otimes_kK$, the corresponding Hopf algebra 
$A = K[\ers\otimes_kK]^*$ is isomorphic as an algebra to the
restricted enveloping algebra, $u(\gl)$, of an abelian Lie algebra $\gl$ 
over $K$ with trivial $p$-mapping.  Given an $A$-module $N$, let 
$\tilde{N}$ denote the module $N$ considered as a $u(\gl)$-module.  
As Yoneda algebras, we necessarily have 
$\Ext_A^*(N,N) \simeq \Ext_{u(\gl)}^*(\tilde{N},\tilde{N})$.  
Hence, an element of $\Ext_A^*(N,N)$ is nilpotent if and only if it is 
nilpotent when considered as an element of
$\Ext_{u(\gl)}^*(\tilde{N},\tilde{N})$.  As nilpotence is detected by 
certain Hopf subalgebras of the form $K[x]/(x^p)$ in the latter case, this 
is also true in the former. So, if $z$ is nilpotent upon restriction to all 
subalgebras $K[x]/(x^p)$, it is nilpotent upon restriction to each 
$\ers\otimes_kK$ and hence, by Corollary \ref{ext-version}, is itself 
nilpotent.
\end{proof}

\begin{rem} Evidently it suffices to take only some of the subalgebras of the 
form $K[x]/(x^p)$, i.e., those which are ``contained in'' an elementary
subgroup scheme of $G$.
\end{rem}


\section{Detecting Projectivity}

As seen in the previous section, an almost immediate corollary of the 
restriction theorem on nilpotent cohomology classes 
(Theorem \ref{restriction}) is a generalization of Chouinard's
theorem (Theorem \ref{Chouinard}) for finite dimensional modules 
(Corollary \ref{proj-cor}). In this section, we show that essentially the 
same proof as for Theorem 6.1 can be used to extend this result to infinite 
dimensional modules, by replacing applications of Proposition
\ref{QV} with applications of Lemma \ref{action}.   

\begin{thm}\label{proj} Let $k$ be a field of characteristic $p > 0$, $G$ be 
a finite group scheme over $k$, and $M$ be a rational $G$-module.  Suppose 
further that the {\rm{(}}infinitesimal{\rm{)}} connected component of the 
identity, $G_0$, of $G$ is unipotent and that all the points of $G$ are 
$k$-rational.  Then $M$ is projective as a $G$-module if and
only if for every field extension $K/k$ and closed subgroup scheme 
$H \subset G\otimes_kK$ with $H \simeq \ers\otimes_kK$ the restriction of 
$M\otimes_kK$ to $H$ is projective.
\end{thm}

\begin{rem} If $M$ is assumed to be finite dimensional, then the proof shows 
that it suffices to take the single field extension $K = \bar{k}$, the 
algebraic closure of $k$.  In this sense, it is a slightly stronger result 
than Corollary \ref{proj-cor} and more comparable to
Proposition 7.6 of \cite{SFB2} (as well as to Chouinard's theorem).
\end{rem}

\begin{proof} 
As previously noted, if $M$ is projective over $G$, then it remains so upon 
restriction to any closed subgroup scheme (cf.~\cite{Jan1}).  Conversely, 
suppose that all restrictions are projective.  The outline  of the proof is 
the same as for the proof of Theorem \ref{restriction} with the details 
modified along the lines of the proof of the Theorem in
\cite{Bend} which was based on the original arguments of Chouinard \cite{Ch}. 

By the assumption on the points of $G$ (see the Introduction), we may write 
$G = G_0\rtimes\pi$ as usual, with $G_0$ assumed to be unipotent.  The first 
step is to reduce 
to the case that $\pi$ is a $p$-group.  Let $\pi_S \subset \pi$ be a 
$p$-Sylow subgroup of $\pi$, and let $G_S$ denote the subgroup 
$G_0\rtimes\pi_S \subset G$.  To show that $M$ is projective, it suffices 
to show that
$\Ext_G^i(M,N) \simeq H^i(G,\Hom_k(M,N)) = 0$ for all $i > 0$ and any 
rational $G$-module $N$.  If $M$ is projective over $G_S$, then 
$H^i(G_S,\Hom_k(M,N)) = 0$ for all $i > 0$, and
since, by Lemma \ref{injection}, the restriction map
$H^*(G,\Hom_k(M,N)) \to H^*(G_S,\Hom_k(M,N))$ is an injection, we also have
$H^i(G,\Hom_k(M,N)) = 0$ for all $i > 0$.  Thus it suffices to assume that 
$\pi$ is a $p$-group and $G$ is in fact unipotent.

Under the assumption that $G$ is unipotent, to show that $M$ is projective 
over $G$, it suffices, by Proposition \ref{uniph1}, to show that 
$H^1(G,M) = 0$. If $L/k$ is any field extension, then 
$H^*(G\otimes_kL,M\otimes_kL) = H^*(G,M)\otimes_kL$ (cf.~\cite{Jan1}).  
Hence, if it can be shown that $H^1(G\otimes_kL,M\otimes_kL) = 0$ for some 
field extension $L/k$, then we also have $H^1(G,M) = 0$.  So it suffices to 
assume that $k$ is algebraically closed.
 
We now proceed by induction on $\dim_kk[G]$ and are trivially done if 
$\dim_kk[G] = 1$ or $G = \ers$, and so assume that the theorem holds for 
all groups $H$ over any field $K/k$ with $\dim_KK[H] < \dim_kk[G]$.

As in the proof of Theorem \ref{restriction}, the next step is to reduce the 
problem to the case that $\pi$ is elementary abelian. If $\pi$ is not 
elementary abelian, by Serre's theorem (Theorem \ref{Serre}), there exists a 
finite product 
$u = \prod u_i = \prod\beta(v_i) \in H^*(\pi,k)$ which is zero in 
$H^*(\pi,k)$ for non-zero $v_i \in H^1(\pi,k)$.  Each $v_i$
can be considered as a non-zero map $\pi \onto \zz$, which can be extended 
to $\phi_i: G \onto \pi \onto \zz$. Let $N_i$ denote the kernel of $\phi_i$.  
Since $\phi_i$ is non-trivial, $\dim_kk[N_i] < \dim_kk[G]$ and so by 
induction $M$ is projective upon restriction to $N_i$ for each $i$.  Hence, 
$H^j(N_i,M) = 0$ for all $j > 0$ and each $i$, and so by 
Lemma \ref{action} the action of $\phi_i^*(z_1) = \psi^*(u_i)$ (where 
$\psi^*: H^*(\pi,k) \to H^*(G,k)$ is induced from the canonical projection
$\psi: G \to \pi$) induces a periodicity isomorphism
$H^{j}(G,M) \to H^{j+2}(G,M)$ for all $j > 0$.  Hence, the product 
$\prod_{i}\phi_i^*(u_i) = \psi^*(u) \in H^*(G,k)$ acts via an isomorphism
on $H^1(G,M)$.  But, $\psi^*(u)$ is zero and so we must have 
$H^1(G,M) = 0$.  Thus $M$ is projective.

From now on we may suppose that $\pi = E_s$ for some $s$ (possibly zero).
In the trivial case $s = 0$, we have $G = G_0$ and the claim is precisely
the Theorem in \cite{Bend}.  In any case, 
the following argument still holds (with Case I being impossible in that 
situation).  As in \ref{restriction}, the rest of the proof consists of
three steps based on $\dim_k\Hom_{Gr/k}(G_0,\ga)^{\pi}$ with arguments like 
the preceding one.

\smallskip
\noindent
{\bf CASE I:} $\dim_k\Hom_{Gr/k}(G_0,\ga)^{\pi} = 0$

Consider the natural projection $\psi: G = G_0\rtimes \pi \to \pi$ and the 
induced map $\psi^*: H^*(\pi,k) \to H^*(G,k)$.  By assumption and 
Lemma \ref{h1}, we have an isomorphism $\psi^*: H^1(\pi,k) \to H^1(G,k)$.  
Consider the map in degree 2, $\psi^*: H^2(\pi,k) \to H^2(G,k)$.  By 
Lemma \ref{h1h2}, either $G \simeq \pi \simeq E_s$ and we are done
or $I \equiv \ker(\psi^*)$ contains some non-zero element of degree 2.  
Again, the map $\psi^*$ preserves the action of the Steenrod algebra and 
hence $I$ is a homogeneous ideal stable under the action of the Steenrod 
operations.  By Proposition \ref{Serre1}, there exists a non-zero product 
$u = \prod u_i = \prod\beta(v_i) \in H^*(\pi,k)$ which lies in $I$.  
In other words, its image under $\psi^*$ is zero in $H^*(G,k)$.  On the
other hand, applying the above argument to $u$, we conclude that $\psi^*(u)$
acts isomorphically on $H^1(G,M)$ and hence $H^1(G,M) = 0$.

\smallskip
\noindent
{\bf CASE II:} $\dim_k\Hom_{Gr/k}(G_0,\ga)^{\pi} = 1$

Just as in \ref{restriction}, there is a homomorphism 
$\eta: G \to \gao\times E_s = \eos$ which is either an isomorphism 
(in which case we're done) or there exists a non-zero element of
degree 2 in $I = \ker\{\eta^*: H^*(\eos,k) \to H^*(G,k)$.  In the latter 
case, by Proposition \ref{Serre3}, there exists a product 
$u = x_1^n\prod u_i \in H^*(\eos,k)$ which lies in $I$.  

Each element $u_i$ may again be identified with $\beta(v_i)$ for some 
non-zero $v_i \in H^1(E_s,k)$.  Just as above, for each $i$, we consider 
the homomorphism $\phi_i: G \onto E_s \onto \zz$ corresponding to $v_i$ 
and conclude that the action of $\eta^*(u_i) = \phi_i^*(z_1) \in H^2(G,k)$ on
$H^*(G,M)$ gives an isomorphism $H^j(G,M) \xto{\sim} H^{j+2}(G,M)$ for all 
$j > 0$.

On the other hand, consider the composite
$$
\phi: G = G_0 \rtimes E_s \xto{\eta} \gao\times E_s \onto \gao,
$$
where the last map is the canonical projection. By Lemma \ref{action}, 
$\phi^*(x_1) = \eta^*(x_1)$  induces a periodicity isomorphism 
$H^j(G,M) \xto{\sim} H^{j+2}(G,M)$ for all $j > 0$.  Hence $\eta^*(u)$
acts isomorphically on $H^1(G,M)$.
Since $\eta^*(u)$ is zero in $H^*(G,k)$ we must
have that $H^1(G,M) = 0$.

\smallskip
\noindent
{\bf CASE III:} $\dim_k\Hom_{Gr/k}(G_0,\ga)^{\pi} > 1$

This case is again divided into two sub-cases based on
$\dim_k\Hom_{Gr/k}(G_0,\gao)^{\pi}$.

\smallskip
{\bf CASE III(a):} $\dim_k\Hom_{Gr/k}(G_0,\gao)^{\pi} = 1$

As in \ref{restriction}, we deduce that there is a homomorphism 
$\eta: G \to \gar\times E_s = \ers$ which is either an isomorphism 
(in which case we're done) or the kernel of 
$\eta^*: H^*(\ers,k) \to H^*(G,k)$ contains a non-zero element of degree 2.  
In the latter case, by Proposition \ref{Serre3}, the kernel contains an 
element of the form $u = x_r^n\prod u_i$ with $\eta^*(x_r) = \phi^*(x_1)$ 
for some non-trivial $\phi: G \to \gao$ and
$\eta^*(u_i) = \phi_i^*(z_1)$ for some non-trivial $\phi_i: G \to \zz$ 
for each $i$.  As above, the actions of $\eta^*(x_r)$ and $\eta^*(u_i)$ 
(for each $i$) induce a periodicity isomorphism $H^j(G,M) \to H^{j+2}(G,M)$ 
for all $j > 0$. Hence, $\eta^*(u)$ acts isomorphically, but since it's zero 
we must have $H^1(G,M) = 0$.

\smallskip
{\bf CASE III(b):} $\dim_k\Hom_{Gr/k}(G_0,\gao)^{\pi} > 1$

Let $\theta_1$ and $\theta_2$ be two linearly independent and non-trivial 
homomorphisms $G_0\to \gao$ which preserve the action of $\pi$.  These can 
then be extended to two linearly independent and non-trivial homomorphisms 
$\phi_1, \phi_2 : G \to \gao$. We now apply the
same argument as in \cite{Bend} to the maps $\phi_1$ and $\phi_2$.  

For any homorphism $\phi: G \to \gao$, let $x_{\phi}$ denote the image
of the canonical generator $x_1 \in H^2(\gao,k)$ under the induced 
map in cohomology.  Once again, the inductive argument as above shows that 
both $x_{\phi_1} = \phi_1^*(x_1)$ and $x_{\phi_2} = \phi_2^*(x_1)$ induce 
periodicity isomorphisms $H^j(G,M) \xto{\sim} H^{j+2}(G,M)$ for all $j > 0$.  
Moreover, for any $c_1, c_2 \in k$, by Corollary 1.5 of \cite{SFB2}, 
$$
c_1x_{\phi_1} + c_2x_{\phi_2} = x_{c_1^{1/p}\phi_1} + x_{c_2^{1/p}\phi_2} 
= x_{c_1^{1/p}\phi_1 + c_2^{1/p}\phi_2} \in H^2(G,k).
$$
If at least one of $c_1$, $c_2$ is non-zero, by linear independence, the 
map $c_1^{1/p}\phi_1 + c_2^{1/p}\phi_2: G \to \gao$ is non-trivial and hence 
by Lemma \ref{action}, the element $c_1x_{\phi_1} + c_2x_{\phi_2}$ also 
induces a periodicity isomorphism $H^1(G,M) \xto{\sim} H^3(G,M)$.  On the 
contrary, since $k$ is by assumption algebraically closed, if $M$ is
finite-dimensional, an eigenvalue argument implies that for some 
$c_1, c_2$ with not both zero, this map will not be an isomorphism unless 
$H^1(G,M) = 0$.  Thus, for finite dimensional modules the proof is complete 
and we see that it is only necessary to extend to the algebraic
closure of $k$.

For infinite dimensional $M$, this eigenvalue argument need not work and 
must be replaced by an infinite dimensional substitute used in \cite{BCR2} 
and which requires a further field extension.  Let $K/k$ be any non-trivial 
algebraically closed field extension.  After base change, the (extended) 
maps $\phi_1, \phi_2: G\otimes_kK \to \gao\otimes_kK$ remain linearly 
independent.  Hence, as the inductive arguments apply equally
well over $K$, we again conclude that for any $c_1, c_2 \in K$ with not 
both zero, the element
$$
c_1x_{\phi_1} + c_2x_{\phi_2} = x_{c_1^{1/p}\phi_1 + c_2^{1/p}\phi_2} 
\in H^2(G\otimes_kK,K)
$$
also induces a periodicity isomorphism 
$$
H^1(G\otimes_kK,M\otimes_kK) = H^1(G,M)\otimes_kK \to
H^3(G\otimes_kK,M\otimes_kK) = H^3(G,M)\otimes_kK.
$$
Applying Lemma 4.1 of \cite{BCR2} (or Lemma 2 of \cite{Bend}) we conclude 
that $H^1(G,M) = 0$ as desired.
\end{proof}

\begin{ques} Can this result be extended to any finite group scheme $G$?
\end{ques}

We end this section by recalling Dade's Lemma \cite[Lemma 11.8]{Dade} which 
says that projectivity (or freeness) of a finite dimensional module for an 
elementary abelian group $E$ can be detected on ``cyclic shifted subgroups'', 
that is, on certain subalgebras of $kE$ of the form $k[x]/(x^p)$.  Indeed, 
this is simply a statement about detecting projectivity (or equivalently
freeness) of modules over a truncated polynomial algebra 
$k[x_1,x_2,\dots,x_n]/(x_1^p,x_2^p,\dots,x_n^p)$.  
However, Dade's lemma as originally stated holds only for finite dimensional 
modules. In \cite{BCR2}, D.~Benson, J.~Carlson, and J.~Rickard extend Dade's 
Lemma to arbitrary modules but at the cost of having to consider field 
extensions.  Indeed, they note that the lemma fails in general without this 
assumption.  We restate it here from a purely algebra perspective.

\begin{prop}[Dade's Lemma -- {\cite[Theorem 5.2]{BCR2}}] Let $k$ be a field 
of characteristic $p > 0$, $n > 0$ be an integer, 
$A_n = k[x_1,x_2,\dots,x_n]/(x_1^p,x_2^p,\dots,x_n^p)$ be a
truncated polynomial algebra over $k$, and $M$ be an $A_n$-module.  
Then $M$ is projective over $A_n$ if and only if for every field extension 
$K/k$, the restriction of $M\otimes_kK$ to the subalgebra 
$K[x]/(x^p) \subset A_n\otimes_kK$ is free for each non-zero 
$x = \sum_{i=1}^nc_ix_i$ with $c_i \in K$.
\end{prop}

\begin{rem} If $M$ is finite dimensional, it suffices to take a single field 
extension $K$ with $K$ algebraically closed.  In general, it suffices to take 
a field extension of transcendence degree at least $n-1$ over the algebraic 
closure $\bar{k}$ of $k$.
\end{rem}

Indeed, it is precisely the necessity of field extensions in this general 
version of Dade's Lemma which makes field extensions necessary in 
Theorem \ref{proj}, whereas they are not necessary for Chouinard's theorem 
(Theorem \ref{Chouinard}).  This is because our generality includes the
case of a restricted Lie algebra.  If we take for example an abelian
Lie algebra with trivial restriction map, then Theorem \ref{proj} becomes 
precisely Dade's Lemma (see the discussion preceding Corollary \ref{xp-nilp}).

Theorem \ref{proj} says that projectivity for $G$ or for the corresponding 
Hopf algebra $k[G]^*$ can be detected on subgroups of the from 
$\ers = \gar\times E_s$.  Denote the corresponding Hopf algebra 
$k[\ers]^* \simeq k[\gar]^*\otimes kE_s$ by $A_{r,s}$.  As an algebra, 
we have already observed that $A_{r,s}$ is simply isomorphic to
$k[x_1,x_2,\dots,x_n]/(x_1^p,x_2^p,\dots,x_n^p)$ where $n = r + s$.  Hence, 
combining Theorem \ref{proj} with Dade's Lemma, we get the following 
detection result.

\begin{cor}\label{xp-proj}  Let $k$ be a field of characteristic $p > 0$, 
$G$ be a finite group scheme over $k$, and $M$ be a rational $G$-module.  
Suppose further that the {\rm{(}}infinitesimal{\rm{)}} connected component 
of the identity, $G_0$, of $G$ is unipotent and that all the points of $G$ 
are $k$-rational.  If for every field extension $K/k$ and subalgebra of the 
form $K[x]/(x^p) \subset K[G\otimes_kK]^*$, $M\otimes_kK$ is free upon
restriction to $K[x]/(x^p)$, then $M$ is projective over $G$.
\end{cor}

\begin{rem} Evidently it suffices to take only some of the subalgebras of 
the form $K[x]/(x^p)$, i.e., those which are ``contained in'' an elementary 
subgroup scheme of $G$.
\end{rem}


\end{document}